\documentclass[12pt]{article}

\usepackage{latexsym}
\usepackage{amsfonts}
\newtheorem{thm}{Theorem}[section]

\newtheorem{lem}[thm]{Lemma}

\title{Deviation Bounds for Wavelet Shrinkage}
\author{ Dawei Hong  \\  
                    \\
  Jean-Camille Birget\thanks{Second author's research
                    supported in part by NSF grant DMS-9970471}\\  \\
	{\footnotesize Dept.\ of Computer Science} \\ 
	{\footnotesize Rutgers University at Camden} \\ 
	{\footnotesize Camden, NJ 08102, USA} \\
        {\footnotesize hong@southwest.msus.edu } \\  
        {\footnotesize birget@camden.rutgers.edu }
       }
\date{}
\begin{document}
\maketitle

\begin{abstract}
We analyse the wavelet shrinkage algorithm of Donoho and Johnstone
in order to assess the quality of the reconstruction of a signal obtained
from noisy samples. We prove {\it deviation bounds} for the maximum
of the squares of the error, and for the average of the squares
of the error, under the assumption that the signal comes from a
H\"older class, and the noise samples are independent, of 0 mean,
and bounded. Our main technique is Talgrand's isoperimetric theorem.
Our bounds refine the known {\it expectations} for the 
average of the squares of the error. 
\end{abstract}

\newpage

%%%%%%%%%%%%%%%%%%%%%%%%%
% Section 1

\section{Introduction}

We address the classical problem of the reconstruction of signal samples
from noisy samples. 
We consider an original {\it signal} of bounded duration \ 
$f$: $t \in [0,1] \to f(t) \in {\mathbb R}$. 
We also have {\it additive noise} $e$: $[0,1] \to {\mathbb R}$. 
Thus, the observed {\it noisy signal} at time $t$ is
$y(t) = f(t) + e(t)$.

We sample the noisy signal at $n$ uniformly spaced instants and we 
denote the sample values by \   
$y_i = f_i + e_i = f(\frac{i}{n}) + e(\frac{i}{n})$ \ 
(for $1 \leq i \leq n$).
Our goal is to recover a good approximation of the original signal 
samples 
$(f_1, \ldots, f_n)$ from the noisy signal samples
$(y_1, \ldots, y_n)$. For this to be possible we need some assumptions
that distinguish the signal from the noise:

\smallskip

\noindent
$\bullet$ \ The original signal $f$ has a certain degree of ``smoothness'', 
i.e., $f$ belongs to a H\"older class
$\Lambda^{\alpha}(M)$ for some $\alpha > 0$ and $M>0$. 

\smallskip

\noindent  
$\bullet$ \ The noise is ``random'', i.e., $(e_1, \ldots, e_n)$ consists
of $n$ independent Borel random variables. 

\medskip

\noindent The {\it H\"older classes} are defined as follows:

\smallskip

\noindent 
For $0 < \alpha \leq 1$, \ $\Lambda^{\alpha}(M)$ \ = \ 
$\{ h \in {\mathbb R}^{[0,1]} \ : \
   (\forall x_1, x_2 \in [0,1]), \
   |h(x_1)-h(x_2)| \leq M |x_1 - x_2|^{\alpha} \}$. 

\smallskip

\noindent
For $1 < \alpha$,  \ $\Lambda^{\alpha}(M) \ = \
\{ h \in {\mathbb R}^{[0,1]} \ : \
   (\forall x \in [0,1]) \ |h'(x)| \leq M, \   
 h^{\lfloor \alpha \rfloor} \ {\mathrm {exists,}} \ {\mathrm {and }}$

$ \ \ \ \ \ \ \ \ \ \ \ \ \ \ \ \ \ \ 
    \ \ \ \ \ \ \ \ \ \ \ \ \ \ \ \ \ \ \ \ \ \ \ \ 
     (\forall x_1, x_2 \in [0,1]) \ 
|h^{\lfloor \alpha \rfloor}(x_1) - h^{\lfloor \alpha \rfloor}(x_2)|
\leq M |x_1 - x_2|^{\alpha - \lfloor \alpha \rfloor} \}$.

\medskip

Let $(\tilde y_1 , \ldots, \tilde y_n )$ be an approximation
of $(f_1, \ldots, f_n)$,  obtained from $(y_1, \ldots, y_n)$. 
Most commonly, the closeness of this approximation is measured by \ 
$\frac{1}{n} \sum_{i=1}^{n} (\tilde y_i - f_i)^2$ \ or by the expectation \ 
${\mathbf E}[\frac{1}{n} \sum_{i=1}^{n} (\tilde y_i - f_i)^2]$ \  
(which makes sense
since the $e_i$, and hence the $\tilde y_i$, are random variables). 

\medskip

The {\bf wavelet shrinkage algorithm} of Donoho and Johnstone \cite{DJ},  
\cite{DJKP} is a very efficient tool for finding good estimates $\tilde y$. 
In outline, the algorithm works as follows: \\ 
(Step 0) \ Choose a wavelet system with $N$ vanishing moments 
($N \geq \alpha$);  
choose a level of coarseness $J_0 \geq 0$ ($J_0$ will depend on $\alpha$),
and consider the multi-resolution chain of Hilbert spaces \ 
$V_{J_0} \subset V_{J_0 +1} \subset \ldots \subset V_j \subset \ldots \ $.  \\ 
(Step 1) \ Apply the {\it Discrete Wavelet Transform} 
(DWT) to the noisy signal samples $(y_1, \ldots, y_n)$, where $n \geq 2^{J_0}$.
This yields the ``empirical wavelet coefficients'' $(\xi_1, \ldots, \xi_n)$.  \\  
(Step 2) \ Fix a ``threshold'' $\lambda_n \ (> 0)$ and apply either ``hard''
or ``soft thresholding'' to $(\xi_1, \ldots, \xi_n)$. 

{\it Hard thresholding} consists of replacing each $\xi_i$ by 0  when
$|\xi_i| \leq  \lambda_n$, and keeping $\xi_i$ unchanged 
when $|\xi_i| > \lambda_n$.

{\it Soft thresholding} consists of transforming each $\xi_i$ as follows: 
$\xi_i$ is replaced by 0 if $|\xi_i| \leq  \lambda_n$; \ if $\xi_i >
\lambda_n$,
$\xi_i$ is replaced by $\xi_i - \lambda_n$; \ if  $\xi_i < -\lambda_n$,
$\xi_i$ 
is replaced by $\xi_i + \lambda_n$.  \\  
(Step 3) \ Apply the inverse DWT to the result of (2). This yields the
estimate 
$(\tilde y_1 , \ldots, \tilde y_n )$.  

\medskip

To what extent does wavelet shrinkage depend on the smoothness conditions 
of the signal $f$ and on the randomness conditions of the noise samples 
$e_i$, and how do the estimators $\tilde y_i$ approximate the original
signal $f$? In \cite{DJ}, \cite{DJKP} it was assumed that the $e_i$ are iid
Gaussian 
variables with distribution {\it N}(0, $\sigma^2$), and the threshold was 
chosen to be $\lambda_n = \sigma \sqrt{2 \frac{\log n}{n}}$. Assuming that 
$f \in \Lambda^{\alpha}(M)$ (the H\"older class) with $\alpha > 0$,  
it is proved in \cite{DJ}, \cite{DJKP} that \  
${\mathbf E}[\frac{1}{n} \sum_{i=1}^{n} (\tilde y_i - f_i)^2] \ < \ 
C \cdot (\frac{1}{n} \log n)^{\frac{2\alpha}{1+2\alpha}}$, \   
where $C$ depends only on $M$ and on the wavelet system used.
It was observed in \cite{DJ}, \cite{DJKP} (the proofs are due to Lepskii 
\cite{Lepskii} and to Brown and Low \cite{BrLow}) that this upper bound is 
optimal over all possible algorithms, if the parameters $\alpha$ and $M$ are
{\it not} known. 
For the optimality of the wavelet shrinkage algorithm it is important that
the threshold be of the form \ $c \cdot \sqrt{\frac{\log n}{n}}$ \  
(where $c$ does not depend on $n$).

\medskip

Since the publication of \cite{DJ}, \cite{DJKP} 
there has been further progress on wavelet shrinkage 
(chapter 6 of \cite{Vi} is an excellent reference up to 1999). 
Most recently, Averkamp and Houdr\'e \cite{AH1}, \cite{AH2} 
expanded the scope of wavelet shrinkage by allowing the noise
samples $e_i$ to have different distributions $F_i$, chosen from a wide
class of distributions. They show in \cite{AH1} (page 32) that the error 
expectation of the wavelet shrinkage algorithm for bounded noise is roughly 
the same as for Gaussian noise, if the parameters $\alpha$ and $M$ of 
the H\"older class of the signal are not known. 
They also discuss various choices of thresholds. 

All the results on wavelet shrinkage in the literature so far evaluate the
quality of the approximation by bounding the expectation  
${\mathbf E}[\frac{1}{n} \sum_{i=1}^{n} (\tilde y_i - f_i)^2]$, to the best
of our knowledge.
In this paper we study deviation bounds (rather than just
the expectation) of \  
$\frac{1}{n} \sum_{i=1}^{n} (\tilde y_i - f_i)^2$ \   
and of \ max$\{ (\tilde y_i - f_i)^2 :  1 \leq i \leq n\}$. 

{\bf Assumptions:} \ We assume that the signal $f$
belongs to a H\"older class $\Lambda^{\alpha}(M)$, and  that the noise 
samples $e_i$ are {\it independent} random variables 
(with possibly different distributions). The only restrictions on the 
distributions are that they are Borel measurable, have {\it compact support} 
(contained in an interval $[- \frac{b}{2}, \frac{b}{2}]$), and zero mean. 
The assumption that the distributions of the noise have bounded support 
is of course equivalent to assuming that the noise $e_i$ has bounded values
($|e_i| \leq \frac{b}{2}$).

The main results of this paper are the following deviation bounds. 

\medskip

\noindent {\bf Theorem.} \ 
{\em 
For the wavelet shrinkage algorithm with threshold 

\smallskip 

\ \ \ \ \ \ \ \ \ \ \ \ \ \ $\lambda_{n, \delta} = C_{\varphi} \, b \, 
(1 + 2 \sqrt{(1+\delta) \ln 2} \, )\sqrt{\frac{\log n}{n}}$  

\smallskip

\noindent
(where $C_{\varphi}$ depends only on the wavelet system) 
we have the following deviation bounds: \\ 
There are $c_1, c_2 > 0$, depending only on $b$, $M$, and $\alpha$, such 
that for all $n \geq n_0$ and all $\delta > 0$,
$$ {\mathbf P} \left({\mathrm max}\{ (\tilde y_i - f_i)^2 : 1 \leq i \leq n\} 
\ \leq \  (c_1 + c_2 \delta) \, 
\left(\frac{\log n}{n}\right)^{\frac{2 \alpha}{1 + 2 \alpha}} \right) 
\ \geq \ 1 - \frac{9}{n^{1 + \delta}}.  $$ 

\smallskip

\noindent
As a consequence,
$$ {\mathbf P} \left(\frac{1}{n} \sum_{i=1}^{n} (\tilde y_i - f_i)^2 
\ \leq \ 
(c_1 + c_2 \delta) \, 
 \left(\frac{\log n}{n} \right)^{\frac{2 \alpha}{1 + 2 \alpha}} \right)
\ \geq \ 1 - \frac{9}{n^{1 + \delta}}.  $$
The minimum number of samples, $n_0$, is $2^9$  when $0 < \alpha \leq 1$; 
when $\alpha > 1$, \\    
$n_0 =  \ (4\alpha+2)^{2\alpha+2} \cdot (\log_2 (4\alpha+2))^2$.
}

\medskip

One notices that  $n_0$ grows very rapidly with $\alpha$, when 
$\alpha > 1$. 
For $\alpha = 2$, we have $n_0 = 1.1 * 10^7$; for 
$\alpha = 3$, \ $n_0 = 3.7 * 10^{10}$, which is impractical. So for large
$\alpha$ our theorem is interesting only from an asymptotic point of view.
On the other hand, in practice usually $\alpha \leq 1$. 

%%%%%%%%%%%%%%%%%%%%%%%%%%%%%%%%%%%%%%%%%%%%%%%%%%%%%%%%%%%%%
%Section 2    

\section{Preliminaries}

\subsection{Wavelets}

We will usually follow the notation of \cite{Daubechies} regarding wavelets,
the only exception being that we reverse the multi-resolution indices. 
Moreover, we only consider real-valued functions with domain $[0,1]$. 
So we have a sequence of real Hilbert spaces \ 
$V_{J_0} \subset V_{J_0 +1} \subset \ldots \subset V_j \subset \ldots \ $,
such that the closure of \ $\bigcup_j V_j$ \ is L$^2[0,1]$.  
We let \ $V_{j+1} = V_j \oplus W_j$ \ (orthogonal complement).
Since we are in the case of compactly supported functions each $V_j$ is a
finite-dimensional real vector space (of dimension $2^j$), with orthonormal 
basis \ $\{ \varphi_{j,k} : 0 \leq k \leq 2^j -1 \}$, derived from a
scaling function $\varphi$.
Let $\psi$ be the wavelet function corresponding to $\varphi$, and let \ 
$\{ \psi_{j,k} : 0 \leq k \leq 2^j -1 \}$ \ be the corresponding orthonormal
basis of $W_j$.  

For any function $g \in \ $ L$^2[0,1]$ we define the piece-wise constant 
function $\overline{g}$: 
$[0,1] \to {\mathbb R}$ \ as follows: \   
$\overline{g}(x) = g(\frac{k}{n}) \  (= g_k)$ \ if
$\frac{k-1}{n} < x \leq \frac{k}{n}$ \ for some $k = 1, \ldots, n$; \
$\overline{g}(x) = 0$ if $x \notin \ ]0,1]$.
The discrete wavelet transform of a vector $(g_1, \dots, g_n)$ can be
obtained by taking the wavelet coefficients of the piecewise constant
function $\overline{g}$. These {\it wavelet coefficients} are:

\smallskip
 
$c_{j,k}^{(g)} = \langle \overline{g}, \varphi_{j,k} \rangle  
               =  \int_0^1 \overline{g}(x) \, \varphi_{j,k}(x) \, dx$, \ and

\smallskip

$d_{j,k}^{(g)} = \langle \overline{g}, \psi_{j,k} \rangle 
               = \int_0^1 \overline{g}(x)\, \psi_{j,k}(x) \, dx$.

\smallskip

\noindent Then for any integer $J \geq J_0$:

\smallskip

$\overline{g}(x) = \sum_{k=0}^{2^J -1} c_{J,k}^{(g)} \varphi_{J,k}(x) + 
        \sum_{j=J}^{+\infty} \sum_{k=0}^{2^j -1} 
                             d_{j,k}^{(g)} \psi_{j,k}(x)$ \ \ \ {\bf a.e.}

\medskip

In this paper we will use two wavelet systems: The Haar wavelets (because 
of their simplicity, especially for programming purposes), and the
interval wavelets with predefined vanishing moments, based on Daubechies
wavelets (Cohen, Daubechies, Jawerth, Vial \cite{CDJV}).

For the {\it Haar wavelets}, the scaling function is \ $\varphi(x) = 1$ 
when $0 < x \leq 1$, and $\varphi(x) = 0$ otherwise. Hence, \ 
$\varphi_{j,k}(x) = 2^{j/2}$ when $k 2^{-j} < x \leq (k+1) 2^{-j}$, and 
$\varphi_{j,k}(x) = 0$ otherwise.
The Haar wavelet function is \ $\psi(x) = 1$ if $0 < x \leq \frac{1}{2}$, \ 
$\psi(x) = -1$ if $\frac{1}{2} < x \leq 1$, and $\psi(x) = 0$ otherwise. 
Hence, \ $\psi_{j,k}(x) = 2^{j/2}$ if 
$k 2^{-j} < x \leq (k+ \frac{1}{2}) 2^{-j}$, \ $\psi_{j,k}(x) = -2^{j/2}$ 
if $(k+ \frac{1}{2}) 2^{-j} < x \leq (k+1) 2^{-j}$, and 
$\psi_{j,k}(x) = 0$ otherwise.

For the {\it interval wavelet system} of \cite{CDJV}, with $N$
{\it vanishing moments}, the scaling function
$\varphi$ and the wavelet function $\psi$ are complicated. But all we 
need to know about them is the following: 

\smallskip

\noindent
$\bullet$ \ A multiresolution of $L^2[0,1]$ is obtained,  with
an orthonormal basis for $V_j$ when $j>J_0$: 

$\{ \varphi_{j,k} : 1 \leq k < 2^j -2N \} \ \cup $
$\{\varphi^{\mathrm{left}}_{j,i}, \varphi^{\mathrm{right}}_{j,i} :$
$0\leq i<N \}.$\\  
Each $\varphi_{j,k}$ has support $[k2^{-j}, (2N-1+k)2^{-j}]$,
each $\varphi^{\mathrm{left}}_{j,i}$ has support $[0, i2^{-j}]$, 
and each $\varphi^{\mathrm{right}}_{j,i}$ has support 
$[1- i2^{-j}, 1]$.  

The decomposition level $J_0$ is chosen so that \ 
$J_0 \geq 1 + \log_2 (2N-1)$. For signals in the H\"older class 
$\Lambda^{\alpha}(M)$ we require the number of vanishing moments to be 
$N \geq \alpha$. 

\smallskip

\noindent
$\bullet$ \ We also have an orthonormal basis for $W_j$, 

$\{ \psi_{j,k} : 1 \leq k < 2^j - 2N \} \ \cup \ $
$ \{\psi^{\mathrm{left}}_{j,i}, \psi^{\mathrm{right}}_{j,i} : 0 \leq i < N
\}$ \\
with the same supports as the corresponding $\varphi$ functions.

\smallskip

\noindent
$\bullet$ \ $\varphi$ and $\psi$ are bounded on $[0,1]$ by a constant
$C >0$, independent of $x$ and $N$: \ \ \  
$\forall x \in [0,1], \ \ $
$|\varphi(x)|, |\psi(x)| \leq C. $

%$C \cdot (1 + |x|)^{-(N+1+\varepsilon)} \ \leq C$.\\   

\medskip

\noindent 
For  $0 \leq k < 2^j - 2N$ (``inside the the interval''), 
$\varphi_{j,k}(x) = 2^{j/2} \varphi(2^j x -k)$. \\  
At the ends of the interval $[0,1]$ we have for $0 \leq i < N$, \ (see
\cite{CDJV}) 

\smallskip

$\varphi^{\mathrm{left}}_{j,i}(x)\ = \ \sum_{h=1}^{2N-1} (-h)^i \varphi(2^j x + h)$.

\smallskip

\noindent
A similar formula holds on the right end of the interval $[0,1]$.
 
\medskip

Assuming that $n$ is a power of 2, $n = 2^J$, we have for the 
function $\overline{y}$, relative to any wavelet system: \ 
$\overline{y}(x) =$
$\sum_{k=0}^{2^J-1} \langle \overline{y}, \varphi_{J,k} \rangle
\varphi_{J,k}(x)$.
Thus for any $J_1$ with $0 \leq J_1 < J$, 
the DWT transforms \ $(y_1, \dots, y_n)$ \ to \  
$\sqrt{n} \, (c_{J_1,0}^{(\overline{y})}, \ldots, $
$c_{J_1,2^J -1}^{(\overline{y})},$
$d_{J_1,0}^{(\overline{y})}, \ldots, d_{J_1,2^J -1}^{(\overline{y})},$ 
$ \ldots, \ldots, $
$ d_{J-1,0}^{(\overline{y})}, \ldots, d_{J-1,2^J -1}^{(\overline{y})})$. 
The DWT is an orthogonal transformation (represented by an orthogonal 
matrix $W$).

We will always assume that $n$ is a power of 2: \ $n = 2^J$.
Throughout this paper, log will refer to log$_2$, and ln will denote
the natural logarithm.

\bigskip

Let us now return to the analysis of a noisy signal $y(t) = f(t) + e(t)$.

\begin{lem}
With respect to the Haar wavelets, the wavelet coefficients of the function
$e$ 
have the following properties: 

\medskip

\noindent
(H1) \ \ For all $j \in [0, 2^J]$ and all $k \in [0, 2^{j-1}-1]$: \ \  
$$c_{j,k}^{(e)} = 2^{-J+j/2} \sum_{i=0}^{2^{J-j}-1} e_{i + 1 + k 2^{J-j}}$$

\bigskip

\noindent
(H2) \ \ For all $j$ and $k$ as in (H1): \ \  
$$d_{j,k}^{(e)} = 2^{-J+j/2} \sum_{i=0}^{2^{J-j-1}-1} 
                         (e_{i + 1 + k 2^{J-j}} 
                          - e_{i + 1 + (k + \frac{1}{2}) 2^{J-j}})$$

\bigskip

\noindent 
For any function $f: [0,1] \to {\mathbb R}$ belonging to
$\Lambda^{(\alpha)}(M)$ with $0< \alpha \leq 1$ we have: 

\medskip

\noindent (H3) \ \ For all $j \in [0, 2^J]$ and all $k \in [0, 2^{j-1}-1]$: 
$$|d_{j,k}^{(f)}| \ < \ M \, 2^{-j(\frac{1}{2} + \alpha)} .  $$ 
\end{lem}

\noindent The proof of this lemma is just a calculation and is given in the 
Appendix.

\medskip

\begin{lem}
With respect to the interval wavelet system \cite{CDJV}, the wavelet coefficients
of the function $e$ have the following properties:

\bigskip

\noindent
(D1) \ \ For all $j \in [0, 2^J]$ and all $k \in [0, 2^{j-1}-1]$: \ \
$$c_{j,k}^{(e)} = 2^{-J+j/2} 
     \sum_{i=0}^{2^{J-j}-1} \alpha_{i,j,k} e_{i + 1 + k 2^{J-j}} $$
\noindent for some numbers $\alpha_{i,j,k}$ that do not depend on the noise
function $e$. 
Moreover, $|\alpha_{i,j,k}| < C_{\varphi}$ \ 
for some constant $C_{\varphi} \geq 1$
depending only on the wavelet system.
 
\bigskip

\noindent
(D2) \ \ For all $j$ and $k$ as in (D1): \ \
$$d_{j,k}^{(e)} = 2^{-J+j/2} \sum_{i=0}^{2^{J-j}-1}
       \beta_{i,j,k} e_{i + 1 + k 2^{J-j}} $$
\noindent for some numbers $\beta_{i,j,k}$ that do not depend on the noise
function $e$. 
Moreover, $|\beta_{i,j,k}| < C_{\varphi}$ \  
where $C_{\varphi} \geq 1 $ depends only on the wavelet system.
\bigskip

\noindent 
Suppose  $f: [0,1] \to {\mathbb R}$ belongs to $\Lambda^{(\alpha)}(M)$ 
with $1 < \alpha$, and suppose the number of vanishing 
moments $N$ of the wavelet system satisfies \ $N \geq \alpha$. 
Then we have:

\medskip

\noindent (D3) \ \ For all $j \in [0, 2^J]$ and all $k \in [0, 2^{j-1}-1]$:
$$|d_{j,k}^{(f)}| \ < \ C_{\varphi} \, M \, 2^{-j(\frac{1}{2} +
\alpha)}$$
where $C_{\varphi} \geq 1$ depends only on the wavelet system.
\end{lem}

\noindent The proof of Lemma 2.2 is just a calculation and is given
in the Appendix.

%%%%%%%%%%%%%%%%%%%%%%%%%%%%%%%%%%%%%%%%%%%%%%%%%%%%%%%%%%%%%%%%%%%%%%%%%%%

\subsection{Talagrand's isoperimetric theorems}

Talagrand's isoperimetric theorems, published in 1995 \cite{Talagrand}, have
had a profound impact on the probabilistic analysis of combinatorial 
optimization methods; Talagrand's theorems often apply quite directly, giving
shorter proofs, often with dramatically better results than previously used 
methods (see \cite{Steele}, chapter 6).  
We will use the following result of \cite{Talagrand}.

Let $(\Omega, \Sigma, \mu_i)$ ($i = 1, \ldots, n$) be Borel 
probability spaces, and
let $\Omega^n$ be the product space with product measure
$P = \mu_1 \times \ldots \times \mu_n$. 
For $A \subseteq \Omega^n$ and 
$\omega = (\omega_1, \ldots, \omega_n) \in \Omega^n$, {\it Talagrand's
`convex' distance} is defined by
$$d_T(\omega,A) = 
\sup \left\{ \ \inf \left\{ \ \sum_{i=1}^n 
\beta_i \cdot I(\omega_i \neq a_i) : 
                              (a_1, \ldots, a_n) \in A  \right\} :   
                     (\beta_1, \ldots, \beta_n) \in {\mathbb R}^n,
                              \sum_{i=1}^n \beta_i^2 = 1    
     \right\}. $$ 
  
\noindent 
Notation:\ $I(\omega_i \neq a_i) = 1$ if $\omega_i \neq a_i$, and 
$I(\omega_i \neq a_i) = 0$ otherwise.

\begin{thm} 
(Talagrand, Theorem 4.1.1 in \cite{Talagrand}): \  For any 
$A \subseteq \Omega^n$ with $P(A) > 0$:
$$\int_{\Omega^n} \exp(\frac{1}{4} d_T(\omega,A)^2) dP(\omega) \ 
\leq \ \frac{1}{P(A)}.$$
As a corollary, for all $t>0$,
$$P(d_T(\omega,A) \geq t) \ \leq \ \frac{1}{P(A)} \cdot
\exp(-\frac{t^2}{4}).$$  
\end{thm}

%%%%%%%%%%%%%%%%%%%%%%%%%%%%%%%%%%%%%%%%%%%%%%%%%%%%%%%%
%Section 3   

\section{Deviation bound for $\frac{1}{n} \sum_{i = 1}^{n}
(f_i - \tilde{y}_i)^2$}

Recall that the input for wavelet shrinkage is $(y_1 , \ldots, y_n )$, 
where $y_i = f_i + e_i$ \ $(i = 1, \ldots, n)$, the $f_i$ are samples 
from the original signal $f$, and the $e_i$ are additive noise.
The $e_i$ are independent Borel random variables.
We assume that the noise is
bounded (with $|e_i| \leq \frac{b}{2}$), so each random variable $e_i$ is
a Borel measurable function $e_i$: 
$\omega_i \in \Omega \mapsto e_i(\omega_i) \in [- \frac{b}{2},
\frac{b}{2}]$. 
Accordingly, we view $(e_1 , \ldots, e_n )$ as a function  \ 
$\omega = (\omega_1, \ldots, \omega_n) \in \Omega^n \ \mapsto \ $
$e(\omega) = (e_1(\omega_1), \ldots, e_n(\omega_n) ) \in [- \frac{b}{2},
\frac{b}{2}]^n$.
(Borel measurability is assumed in order to apply Talagrand's theorem.) 
To simplify the notation we often write $e_i(\omega)$ for $e_i(\omega_i)$.  
 
\medskip

\noindent 
We shall first define a subset $A$ of $\Omega^n$ and then show that
\begin{itemize}
\item 
$P(A) > \frac{1}{9}$ if $n$ is large enough, and
\item wavelet shrinkage satisfies our deviation bounds when the noise 
samples are in $A$.
\end{itemize}
Then for any $\delta > 0$ we define a subset 
$B_{\delta} \subseteq \Omega^n$ such that
\begin{itemize}
\item 
for any $\omega \in \Omega^n$, if  Talagrand's distance satisfies
$d_T (\omega, A) \leq 2 \sqrt{(1 + \delta) \ln{n}}$ \ 
then $\omega \in B_{\delta}$;
\item wavelet shrinkage satisfies our deviation bounds when the noise 
samples are in $B_{\delta}$.
\end{itemize}
Finally, by applying Talagrand's theorem we obtain our results.

%%%%%%%%%

\subsection{The subset $A$}

Recall that we assume $n = 2^J$. 
For any $\omega \in \Omega^n$ we decompose the noise sample sequence 
$e(\omega)$ into blocks of length $J$, as follows: 

$e(\omega) = ( \ldots, \ \ldots, $
              $e_{kJ + 1}(\omega), \ldots, e_{(k + 1)J}(\omega), $
              $\ldots, \ \ldots)$ \\  
where $k = 0, \ \ldots, \ \frac{1}{J} 2^J - 1$. \,  
Here, for simplicity we regard $\frac{1}{J} 2^J = 2^{J - \log J}$ as an
integer (i.e., we assume that $J$ is a power of 2).

\medskip

\noindent 
For the Haar wavelets we define the subset $A \subset \Omega^n$ as follows: 
$$A \ = \ \{\omega \in \Omega^n : ( \forall \ell \in [-1,  J - \log J])
(\forall k \in [0,  2^{J - \log J - \ell} - 1]),  $$ 
$$\left| \sum_{i = 0}^{J2^{\ell-1} - 1} e_{k 2^{\ell} J + i + 1}(\omega) \right| 
\leq  
b J 2^{\ell/2} \sqrt{2^{-1}\ln 2}  \  \}. $$ 

%AAA
\noindent
For the interval wavelet system we define
$$A \ = \ \{\omega \in \Omega^n : ( \forall \ell \in [-1,  J - \log J])
(\forall k \in [0,  2^{J - \log J - \ell} - 1]),  $$
$$\left| \sum_{i = 0}^{J2^{\ell} - 1} e_{k 2^{\ell} J + i + 1}(\omega) \cdot
   \alpha_{i, J-\log J -\ell, k} \right|
\leq
b J 2^{\ell/2} \sqrt{2^{-1}\ln 2}  $$ 
$$ {\mathrm{and}}  \ \ \ \ \  
\left| \sum_{i = 0}^{J2^{\ell} - 1} e_{k 2^{\ell} J + i + 1}(\omega) \cdot
   \beta_{i, J-\log J -\ell, k} \right|
\leq
b J 2^{\ell/2} \sqrt{2^{-1}\ln 2} \  \}. $$

\medskip

\noindent We need a classical result from probability theory.

\begin{thm} 
\label{Hoeffding}
(Hoeffding's inequality)
Let $X_1, \ldots, X_m$ be independent random variables with 
$b_1 \leq X_i \leq b_2$ \  $(i = 1, \ldots, m)$. Then for all $t > 0$,
\[ P \left(\left| \sum_{i = 1}^{m} (X_i - E[X_i]) \right| \leq t \right) 
\geq 1 - \exp \left(- \frac{2t^2}{m(b_2 - b_1)^2} \right). \] 
\end{thm}

\medskip

\begin{lem}
\label{lem1_for_A}
For all $n > 1$, \ $P(A) \geq 1 - \frac{4}{\log n} + \frac{1}{n}$ \ for the
Haar wavelets, and  \ $P(A) \geq 1 - \frac{8}{\log n} + \frac{2}{n}$ \ for
the interval wavelet system.

In either case, if
$n \geq 256$ then $P(A) \geq \frac{1}{128}$.
If $n \geq 2^9$ then $P(A) > \frac{1}{9}$. 
Moreover, $P(A)$ tends to 1 when $n \to \infty$.   
\end{lem}

\noindent
{\bf Proof:} We first give the proof for the Haar wavelets. 
For any $\ell \in [-1, J - \log J]$ and 
$k \in [0, 2^{J - \log J - \ell} - 1]$ the noise samples 
$e_{k2^{\ell} J + 1}, \ldots, e_{(k + 1)2^{\ell} J}$
are independent random variables, each with values in 
$[- \frac{b}{2}, \frac{b}{2}]$. 
So Hoeffding's inequality applies, and since $E[e_i]$ = 0 for all $i$,
we obtain for all $t > 0$,

$$P \left( \left| 
\sum_{i = 0}^{2^{\ell-1} J - 1} e_{k2^{\ell} J + i + 1}
\right| 
\leq t \right)
\ \geq \ 
1 - \exp \left( - \frac{2t^2}{2^{\ell} J b^2} \right). $$

\noindent
Letting \  
$t =  b 2^{\ell /2} J \sqrt{2^{-1} \ln 2}$ \   we obtain 

\begin{equation}
\label{use_Hoeffding2}
P \left( \left| \sum_{i = 0}^{2^{\ell-1} J - 1} e_{k2^{\ell} J + i + 1}
\right| \leq
b 2^{\ell /2} J \sqrt{2^{-1} \ln 2} \right)
\geq 1 - \frac{1}{n}.
\end{equation}

\noindent
For $\ell \in [-1, J - \log J]$ and 
$k \in [0,  2^{J - \log J - \ell} - 1]$, let 
$$ A_{\ell,k} = 
\left\{\omega \in \Omega^n~:~  
\left| \sum_{i = 0}^{2^{\ell-1} J - 1} e_{k 2^{\ell} J + i + 1}(\omega)
\right| \leq 
 b 2^{\ell/2} J \sqrt{2^{-1} \ln 2} \right\}$$
and let \ $A_{\ell} = \bigcap_{k=0}^{2^{J - \log J - \ell} - 1} A_{\ell,k}$.

\smallskip

\noindent
Then by (\ref{use_Hoeffding2}), \   
$P(A_{\ell,k}) \geq 1 - \frac{1}{n}.$  

\smallskip

\noindent For the complements of these sets we have \  
$\bar{A_{\ell}}= \bigcup_{k=0}^{2^{J- \log J -\ell} - 1} \bar{A}_{\ell,k}$

\smallskip

\noindent hence \ 
$P(\bar{A_{\ell}}) \leq  \sum_{k=0}^{2^{J - \log J - \ell} - 1}
\frac{1}{n}$.  

\smallskip

\noindent  
Since $n = 2^J$ we obtain \ 
$P(\bar{A_{\ell}}) \leq \frac{2^{-\ell}}{\log n}$. 

\smallskip

\noindent
Since \  $A = \bigcap_{\ell = -1}^{J - \log J} A_{\ell}$ \  we have

\smallskip

$ P(A) \geq \ 1 - \sum_{\ell = -1}^{J - \log J} P( \bar{A_{\ell}}) \ \geq \ $
$1 - \sum_{\ell=-1}^{J - \log J} \frac{2^{-\ell}}{\log n}.$ 

\smallskip

\noindent  Hence,  \ $P(A) \geq 1 - \frac{4}{\log n} +  \frac{1}{n}.$ \ \ 
This proves the Lemma for the Haar case.

\bigskip

\noindent For the interval wavelet system we let 
$$A^{\alpha} \ = \ 
\{ \omega \in \Omega^n \ : \ ( \forall \ell \in [-1,  J - \log J])
(\forall k \in [0,  2^{J - \log J - \ell} - 1]),  $$
$$ \ \ \ \ \ \ \ \ \ \ \ \ \ \ 
   \left| \sum_{i = 0}^{J2^{\ell} - 1} e_{k 2^{\ell} J + i + 1}(\omega)
\cdot
   \alpha_{i, J-\log J -\ell, k} \right|
\leq
b J 2^{\ell/2} \sqrt{2^{-1}\ln 2} \ \}, $$
and  $$A^{\beta} \ = \ 
\{ \omega \in \Omega^n \ : \ ( \forall \ell \in [-1,  J - \log J])
(\forall k \in [0,  2^{J - \log J - \ell} - 1]),  $$
$$ \ \ \ \ \ \ \ \ \ \ \ \ \ \ 
  \left| \sum_{i = 0}^{J2^{\ell} - 1} e_{k 2^{\ell} J + i + 1}(\omega) \cdot
   \beta_{i, J-\log J -\ell, k} \right|
\leq
b J 2^{\ell/2} \sqrt{2^{-1}\ln 2} \  \}. $$
Then \  $A = A^{\alpha} \cap A^{\beta}$.

\smallskip

\noindent We also let 
$$A^{\alpha}_{\ell,k} = \{ \omega \in \Omega^n \ : \   
\left| \sum_{i = 0}^{J2^{\ell} - 1} e_{k 2^{\ell} J + i + 1}(\omega) \cdot
   \alpha_{i, J-\log J -\ell, k} \right|
\leq
b J 2^{\ell/2} \sqrt{2^{-1}\ln 2} \ \}, $$
and $$A^{\beta}_{\ell,k} = \{ \omega \in \Omega^n \ : \ \left| \sum_{i =
0}^{J2^{\ell} - 1} e_{k
2^{\ell} J + i + 1}(\omega) \cdot
   \beta_{i, J-\log J -\ell, k} \right|
\leq
b J 2^{\ell/2} \sqrt{2^{-1}\ln 2} \  \}. $$
Moreover, we let \ $A^{\alpha}_{\ell} = \bigcap_k A^{\alpha}_{\ell,k}$ \ and
\ 
$A^{\beta}_{\ell} = \bigcap_k A^{\beta}_{\ell,k}$.
Then $A_{\ell} = A^{\alpha}_{\ell} \cap A^{\beta}_{\ell}$, hence 
$\bar{A_{\ell}} = \bar{A^{\alpha}_{\ell}} \cup \bar{A^{\beta}_{\ell}}$.   

\smallskip

\noindent By the same proof as for Haar wavelets above: \ 
$P(\bar{A^{\alpha}_{\ell}})$ and 
$P(\bar{A^{\beta}_{\ell}}) \leq \frac{2^{-\ell}}{\log n}$.   

\smallskip

\noindent 
Hence, \ \  $P(\bar{A_{\ell}}) \leq \frac{2^{-\ell+1}}{\log n}$.  

\smallskip

\noindent   
Since \ $A = \bigcap_{\ell = -1}^{J - \log J} A_{\ell}$ \  we obtain by a
similar calculation as in the Haar case: \ 

\smallskip

$P(A) \geq 1 - \frac{8}{\log n} +  \frac{2}{n}$.
\ \ \ \ \ \ \ \ \       $\Box$ 

\bigskip

\begin{lem}
\label{lem2_for_A}
For all $\omega \in A$, all $j \in \ ]J_0, J[$, and all
$k \in [0, 2^{j} - 1]$, we have (for some constant 
$C_{\varphi} \geq 1$,
depending only on the wavelet system):

$$ |d_{j, k}^{(e(\omega))}| \leq 
 b \, C_{\varphi} \, \sqrt{ \frac{\log n}{n} } $$
and for all $k \in [0, 2^{J_0} - 1]$,
$$ |c_{J_0, k}^{(e(\omega))}| \leq 
 b \, C_{\varphi} \, \sqrt{\frac{\log n}{n} } $$
\end{lem}
 
\noindent
{\bf Proof:} We consider two cases for $j$. \\  
Case 1: \ $J_0 \leq j \leq J - \log J + 1$.  \\   
We write $j$ as
$J - \log J - \ell$, where  $-1 \leq \ell \leq J - \log J - J_0$. 
Let us first consider Haar wavelets. By (H2) (in Lemma 2.1) we have
$$ d_{j, k}^{(e(\omega))} = 2^{-J+j/2}
\left(
\sum_{i = 0}^{2^{\ell - 1}J - 1}
  e_{k2^{\ell}J + i + 1}(\omega) - 
\sum_{i = 0}^{2^{\ell - 1}J - 1} 
  e_{(k + 1/2)2^{\ell}J + i + 1}(\omega) 
\right). $$ 
Since $\omega \in A$ we can apply the defining property of $A$ to  

\smallskip

$\left| \sum_{i=0}^{J2^{\ell-1}-1} e_{i + 1 + k2^{\ell}J} \right| \ = \  $
$\left| \sum_{i=0}^{J2^{\ell-1}-1} e_{i + 1 + 2k2^{\ell-1}J} \right|$.  

\smallskip

\noindent  Since $2k$ is in the correct range \ 
$[0, 2^{j+1}-2] = [0, \frac{1}{J}2^{J-(\ell-1)}-2]$, \   we have  

\smallskip

$\left| \sum_{i=0}^{J2^{\ell-1}-1} e_{i + 1 + k2^{\ell}J} \right| \leq $
$bJ2^{(\ell-1)/2} \sqrt{2^{-1}\ln 2}$. 

\smallskip

\noindent Similarly, 

\smallskip

$\left| \sum_{i=0}^{J2^{\ell-1}-1} e_{i + 1 + (k+\frac{1}{2})2^{\ell}J}
\right| = $
$\left| \sum_{i=0}^{J2^{\ell-1}-1} e_{i + 1 + (2k+1)2^{\ell-1}J} \right| $
$\leq bJ2^{(\ell-1)/2} \sqrt{2^{-1}\ln 2}$ \, ;

\smallskip

\noindent we used the defining property of $A$, since the range of $2k+1$ is

\smallskip

$[0, 2^{j+1}-2+1] = [0, \frac{1}{J}2^{J-(\ell-1)}-1]$.    

\smallskip

\noindent By combining these two bounds we obtain   

\smallskip

$|d_{j, k}^{(e(\omega))}|  \ \leq \    $
$2^{-J+j/2} \cdot 2 \cdot bJ2^{(\ell-1)/2} \sqrt{2^{-1}\ln 2} \ < \  $
$b \sqrt{\ln 2} \sqrt{\frac{\log n}{n}} \ \leq \ $   
$b \sqrt{\frac{\log n}{n}}$.

\bigskip

Let us now consider case 1 for the interval wavelet system. By (D2) in
Lemma 2.2, 

\smallskip

\ \ \ \ \ \ $ d_{j, k}^{(e(\omega))} = 2^{-J+j/2} \cdot 
\sum_{i = 0}^{2^{\ell}J - 1}
  e_{k2^{\ell}J + i + 1}(\omega) \cdot \beta_{i,j,k}. $

\smallskip

\noindent
Since $\omega \in A$, 

$|d_{j, k}^{(e(\omega))}| \ \leq \ 2^{-J+j/2} \cdot bJ2^{(\ell-1)/2}
\sqrt{2^{-1}\ln 2} $
$\ = \ b2^{(-J + \log J)/2} \sqrt{2^{-1}\ln 2}$
$\ = \  b \sqrt{\frac{\log n}{n}} \sqrt{2^{-1}\ln 2} $

$\ \ \ \ \ \ \leq \ b \sqrt{\frac{\log n}{n}}.$

\bigskip

\noindent Case 2: \  $J - \log J + 2 \leq j < J$. \\   
For the Haar wavelets we use the boundedness of the noise, 
$|e_i - e_j| \leq b$. 
Hence, by (H2), 

\smallskip

$ |d_{j ,k}^{(e(\omega))}| \leq 2^{-J+j/2} b (J2^{\ell-1} -1)  \ \leq \ $  
$ b \sqrt{\frac{\log n}{n}}.$ 

\medskip

\noindent For the interval wavelet system, (D2) yields

\smallskip

$|d_{j ,k}^{(e(\omega))}| \leq \ 2^{-J+j/2}$
$\sum_{i=0}^{2^{J-j}-1} |e_{k2^{\ell}J + i + 1}(\omega)| \cdot
|\beta_{i,j,k}|$ 
$= 2^{-J+j/2} \, 2^{J-j} \, \frac{b}{2} \, C_{\varphi} $ 

$\leq  \ \frac{b}{2} \, C_{\varphi} \, 2^{-j/2} \ \leq $
$\ b \, C_{\varphi} \sqrt{\frac{\log n}{n}}$

\smallskip

\noindent by using $j \geq J - \log J +2$ for the last inequality. 

\medskip

\noindent
By an argument similar to the above we obtain the bound for 
$|c_{J_0,k}^{(e(\omega))}|$. \ \ \ \ \  $\Box$

\bigskip

\bigskip

To implement wavelet shrinkage we need two parameters: A decomposition level
$J_0$ and a threshold $\lambda_{n,\delta}$. We define 
$$J_1 = \lceil \frac{1}{1 + 2\alpha}
(J - \log J) \rceil$$
and we choose $J_0$ so that \ \  $J_0 \leq J_1.  $

For the Haar wavelets (when $0 < \alpha \leq 1$) we can simply pick $J_0 =0$, 
but for the interval wavelet system (when $1 < \alpha$ and we have 
$N = \lceil \alpha \rceil$ vanishing moments), we also require (see
\cite{CDJV}) that \  $J_0 \geq 1 + \log(2N-1)$.
When $\alpha > 1$ we choose 
$$J_0 = 1 + \lceil \log(2 \, \lceil \alpha \rceil -1) \rceil$$
Thus, for $J_0$ to exist (when $\alpha > 1$) we need $n = 2^J$ to be such that \  
$1 + \log(2 \lceil \alpha \rceil -1) \leq J_1$. \ 
A sufficient condition for this is that \ 
$J - \log J  \geq (1+ \log(2\alpha +1) ) \, (1+ 2\alpha) $, \\   
or equivalently, \ \ \ \ 
$\frac{n}{\log n} \geq (4\alpha +2)^{2\alpha +1} $.

\smallskip

\noindent By using the fact that \ $\frac{n}{\log n}$ \ is an increasing
function of $n$, 
and that the relation \  $\frac{y}{\log y} \geq x$ \ is implied by \ 
$y \geq \ x \cdot \log x \cdot \log \log x$, \ we have the following
sufficient condition
on $n$: 

\medskip

When $\alpha > 1$ we assume that 
$$n \ \geq \ 
(4\alpha+2)^{2\alpha+2} \cdot (\log(4\alpha+2))^2$$

\medskip

\noindent We use the threshold
$$\lambda_{n,\delta} = C_{\varphi} \, b \, 
      \left(1 + 2 \sqrt{(1+\delta) \ln 2} \right) 
                                     \sqrt{ \frac{\log n}{n} }$$

\medskip

The first step of the wavelet shrinkage algorithm is DWT, which maps
$(y_1 , \ldots, y_n )$ to  \
$\sqrt{n} \, (c_{J_0 , 0}^{(y)}, \ldots, c_{J_0 , 2^{J_0} - 1}^{(y)},
d_{J_0, 0}^{(y)}, \ldots, d_{J_0, 2^{J_0} - 1}^{(y)}, \ldots, \ldots,
d_{J - 1, 0}^{(y)}, \ldots, d_{J - 1, 2^{J - 1} - 1}^{(y)})$, \
where $n = 2^J$. \\
Since $y_i = f_i + e_i$ and the DWT is linear we have
$$c_{J_0 , k}^{(y)} = c_{J_0 , k}^{(f)} + c_{J_0 , k}^{(e)}, ~~~~~0 \leq k <
 2^{J_0},$$
and
$$d_{j, k}^{(y)} = d_{j, k}^{(f)} + d_{j, k}^{(e)}, ~~~~~ J_0  \leq j <
 J,~~~0 \leq k < 2^j,$$
where
$c_{J_0 , k}^{(f)}$, $d_{j, k}^{(f)}$ and $c_{J_0 , k}^{(e)}$, $d_{j,
k}^{(e)}$
are the wavelet coefficients for $(f_1 , \ldots, f_n )$
and $(e_1 , \ldots, e_n )$, respectively.

\medskip

The second step of wavelet shrinkage is thresholding. We shall
prove our result for soft thresholding. But in our proofs it will be
easy to see that our results will hold for hard thresholding too.
For soft thresholding, we have
\[ \tilde d_{j ,k} =  \left\{   \begin{array}{ll}
 d_{j ,k}^{(y)} - \lambda_{n,\delta} & \mbox{if \ $d_{j ,k}^{(y)} >
\lambda_{n,\delta}$} \\
 0                    & \mbox{if \ $|d_{j ,k}^{(y)}| \leq
\lambda_{n,\delta}$} \\
 d_{j ,k}^{(y)} + \lambda_{n,\delta} & \mbox{if \ $d_{j ,k}^{(y)} <
-\lambda_{n,\delta}$}
                              \end{array}
                     \right.                 \]

\medskip

The last step of wavelet shrinkage is the inverse of DWT which yields
$\tilde y = (\tilde y_1 , \ldots, \tilde y_n )$.
If we let
\begin{equation}
\label{def_of_tildey}
\tilde y(x)= \sum_{k=0}^{2^{J_0}-1} c_{J_0,k}^{(y)} \, \varphi_{J_0,k}(x)
+ \sum_{j=J_0}^{J-1} \sum_{k=0}^{2^{j-1}} \tilde d_{j,k} \, \psi_{j,k}(x),
\end{equation}
then we obtain \ $\tilde y_i = \tilde y (\frac{i}{n})$ \ for 
$i = 1, \ldots,n$.

\medskip

\subsection{Application of Talagrand's theorem}

Let $W$ be the orthogonal matrix that represents the DWT. Let 
$A \subseteq \Omega^n$ be as above. 
For any $\delta > 0$ we define the following subset of $\Omega^n$:
$$B_{\delta} = \left\{\omega' \in \Omega^n : \   
   (\forall \ell \in [1,n]), \ \   
   \inf_{\omega \in A} \left| \sum_{i=1}^n W_{\ell,i}(e_i(\omega') -
e_i(\omega))\right|
    < \ 2 b \sqrt{(1+\delta) \ln n} \ \right\}. $$

\begin{lem} 
For all $\omega' \in B_{\delta}$ and all $k \in [0, 2^{J_0}-1]$: \ \   
$|c^{(e(\omega'))}_{J_0,k}| \leq  \lambda_{n,\delta}$. 

\medskip

\noindent For all $j \in [J_0, J-1]$ and $k \in [0, 2^j -1]$: \ \ 
$|d^{(e(\omega'))}_{j,k}| \leq \lambda_{n,\delta}$. 
\end{lem}

\noindent
{\bf Proof:} By the definition of $B_{\delta}$, for every 
$\omega' \in B_{\delta}$ there exists $\omega \in A$ such that 
$$ \sqrt{n} \ |c^{(e(\omega))}_{J_0,k} - c^{(e(\omega'))}_{J_0,k}| \leq 
b 2 \sqrt{(1 + \delta) \ln n} $$
and 
$$ \sqrt{n} \ |d^{(e(\omega))}_{j,k} - d^{(e(\omega'))}_{j,k}| \leq 
b 2 \sqrt{(1 + \delta) \ln n} $$ 
The Lemma then follows from Lemma \ref{lem2_for_A}.  \ \ \ \ $\Box$

\medskip

For the following theorem we use the threshold $\lambda_{n,\delta}$ as
above;
we let $n_0 = 2^9$ when $0 < \alpha \leq 1$, and  \  
$n_0 = \ (4\alpha+2)^{2\alpha+2} \cdot (\log(4\alpha+2))^2$ \  
when $\alpha > 1$.

\begin{lem}
When $n \geq n_0$, \ \ 
$P(B_{\delta}) > 1 - \frac{9}{n^{1+ \delta}}$. 
\end{lem}
{\bf Proof:} We first prove that 
$$\{ \omega' \in \Omega^n \ : \ 
d_T(\omega',A) < 2\sqrt{(1+\delta) \ln n}\, \} \ \subseteq \ B_{\delta}.$$
Recall the definition \\  
$d_T(\omega',A) = $ \\  
$\sup \{ \ \inf \{ \sum_{i=1}^n \beta_i \cdot I(\omega'_i \neq \omega_i) :$
$(\omega_1, \ldots, \omega_n) \in A \} : $
$(\beta_1,\ldots,\beta_n)\in {\mathbb R}^n,$
                 $ \sum_{i=1}^n \beta_i^2 = 1 \}. $

\smallskip

\noindent
We will choose the following $n$ vectors for
$\beta = (\beta_1, \ldots, \beta_n)$ in the above formula:

\smallskip

$(|W_{1,\ell}|, \ldots, |W_{n,\ell}|)$,  for $\ell = 1, \ldots, n$.

\smallskip

\noindent
Since $W$ is orthogonal all  its row vectors have unit length.
For all $\omega' \in \Omega^n$, $\omega = (\omega_1, \ldots, \omega_n) \in A$,
and $1 \leq \ell \leq n$, we have:

\smallskip

$| \sum_{i=1}^{n} W_{i,\ell} (e_i(\omega') - e_i(\omega)) | $

$\leq \ \ b \sum_{i=1}^{n} |W_{i,\ell}| \cdot I(e_i(\omega') \neq e_i(\omega))$

$\leq \ b \sum_{i=1}^{n} |W_{i,\ell}| \cdot I(\omega' \neq \omega)$.

\smallskip

\noindent
(The last inequality follows from the fact that \
$I(e_i(\omega') \neq e_i(\omega)) \leq I(\omega' \neq \omega)$, \
because $e_i(\omega') \neq e_i(\omega)$ implies $\omega' \neq \omega$.)

\smallskip

\noindent
Hence, for all $\omega' \in \Omega^n$  and $1 \leq \ell \leq n$,

\smallskip

$\inf \{ |\sum_{i=1}^{n} W_{i,\ell} (e_i(\omega') - e_i(\omega))| : $
$\omega \in A\} $

$ \leq \ \inf \{ \sum_{i=1}^{n} |W_{i,\ell}| \cdot I(\omega' \neq \omega) \, b$
$ : \omega \in A \}$

$= \ b \ \inf \{ \sum_{i=1}^{n} |W_{i,\ell}| \cdot I(\omega' \neq \omega)$
$ : \omega \in A \}$.

\smallskip

\noindent
Therefore, if  \ $d_T(\omega',A) \leq 2 \sqrt{(1+ \delta) \ln n }$ \ \ then
for all $1 \leq \ell \leq n$,

\smallskip

$\inf \{|\sum_{i=1}^{n} W_{i,\ell} (e_i(\omega')-e_i(\omega))| : $
$\omega \in A \} \ $
$\leq \ b \, 2 \sqrt{(1+ \delta) \ln n }.$

\smallskip

\noindent This means that $\omega' \in B_{\delta}$, and this proves that

$\{ \omega' \in \Omega^n \ : \ $
$d_T(\omega',A) < 2\sqrt{(1+\delta) \ln n}\, \} \ \subseteq \ B_{\delta}.$

\smallskip

\noindent Hence, $P(B_{\delta}) \ \geq \  $
$ P(\{ \omega' \in \Omega^n : $
$d_T(\omega',A) < 2\sqrt{(1+ \delta)\ln n} \, \})$.

\smallskip

\noindent
By Talagrand's theorem this is \ 
$\geq \ 1 - \exp(- (1+ \delta) \ln 2) \cdot \frac{1}{P(A)} \ $
$>  \  1 - \frac{9}{n^{1+ \delta}}$.
\ \ \ \ \ \ $\Box$

\begin{lem}
\label{tilda}
For all \ $\omega' \in B_{\delta}$ we have:

\smallskip

\noindent (1) \ \ \ \
When $J_1 \leq j < J$,\ $0 \leq k < 2^j$ , \ \ \ \ \ \ 
$|\tilde d_{j, k}(\omega') - d_{j, k}^{(f)}| \ \leq \ |d_{j,k}^{(f)}| \ $
$\leq \ C_{\varphi} M \cdot 2^{-j(\frac{1}{2} + \alpha)}$. 

\smallskip

\noindent (2) \ \ \ \ 
When $J_0 \leq j < J_1$, \ $0 \leq k < 2^j$,  \ \ \ \ \ \  
$|\tilde d_{j, k}(\omega') - d_{j, k}^{(f)}| \ \leq \ 2\lambda_{n,\delta}$.
\end{lem}
{\bf Proof:}  \ To prove (1), we note first that by (H3), (D3) we have \ 
$|d_{j ,k}^{(f)} | \leq C_{\varphi} M \, 2^{-j(1/2 + \alpha)}$.

\smallskip

\noindent To prove the inequality \ 
$|d_{j, k}^{(f)} - \tilde d_{j, k}| \leq |d_{j, k}^{(f)}|$ \ 
one considers six cases, according to
the possible relative positions of 0, $d_{j,k}^{(f)}$, and $\tilde d_{j,k}$.
If $0 \leq \tilde d_{j, k} \leq d_{j, k}^{(f)}$, or if
$d_{j,k}^{(f)} \leq \tilde d_{j,k} \leq 0$, the inequality is obvious from
the order picture. The other four cases are not possible, since they would
imply that $|d_{j, k}^{(e(\omega))}| > \lambda_{n,\delta}$, contradicting
what we saw a little earlier. This proves (1).

\smallskip

For the proof of (2) we consider two cases. If $\tilde d_{j, k} = 0$,
$|d_{j, k}^{(y)}| \leq \lambda_{n,\delta}$, hence
$|d_{j, k}^{(f)} - \tilde d_{j, k} | = |d_{j, k}^{(f)}| = $
$|d_{j, k}^{(y)} - d_{j, k}^{(e)}| \leq |d_{j, k}^{(y)}| + |d_{j, k}^{(e)}|$
$\leq \lambda_{n,\delta} + \lambda_{n,\delta}$.
In the second case, $|d_{j, k}^{(y)}| > \lambda_{n,\delta}$, and
$|d_{j, k}^{(f)} - \tilde d_{j, k}| = |d_{j, k}^{(e)} - \lambda_{n,\delta}|$
$\leq \lambda_{n,\delta} + \lambda_{n,\delta}$.
This proves the inequality.     \ \ \ \ \ \ $\Box$

\medskip

\begin{thm}
{\bf (Deviation bound for max square error)} \  
For wavelet shrinkage with threshold $\lambda_{n,\delta}$
we have for all $n \geq n_0$:
$$P\left( \max_{0\leq i \leq n} (f_i - \tilde{y_i})^2 \ \leq \  
(c_1 + c_2 \, \delta)
\left(\frac{\log n}{n}\right)^{\frac{2 \alpha}{1 + 2 \alpha}} \right)
\ \geq \  1 - \frac{9}{n^{1 + \delta}} $$
where $c_1$ and $c_2$ depend only on $b$, $M$, and $\alpha$.

\smallskip

\noindent As a consequence {\bf (deviation bound for mean square error)},
$$P\left(\frac{1}{n} \sum_{i=0}^{n} (f_i - \tilde{y_i})^2 \leq
(c_1 + c_2 \, \delta)
\left(\frac{\log n}{n}\right)^{\frac{2 \alpha}{1 + 2 \alpha}} \right)
\ \geq \  1 - \frac{9}{n^{1 + \delta}} $$
\end{thm}
{\bf Proof:} \ At the beginning of subsection 2.1 we defined the function
$\overline{f}$, and its wavelet coefficients. We have  
\[\overline{f}(x) =
\sum_{k = 0}^{2^{J_0} - 1} c_{J_0 , k}^{(f)} \varphi_{J_0 , k}(x) +
\sum_{j = J_0}^{J_1 - 1} \sum_{k = 0}^{2^j - 1}
d_{j, k}^{(f)} \psi_{j, k}(x) +
\sum_{j = J_1}^{J - 1} \sum_{k = 0}^{2^j - 1}
d_{j, k}^{(f)} \psi_{j, k}(x),
\]
and \ $f_i = \overline{f}(\frac{i}{n})$ \ for $1 \leq i \leq n$.

In connection with the thresholding of $y$ we define the function 
\[\tilde y(x) =
\sum_{k = 0}^{2^{J_0} - 1} c_{J_0 , k}^{(y)} \varphi_{J_0 , k}(x) +
\sum_{j = J_0}^{J_1 - 1} \sum_{k = 0}^{2^j - 1}
\tilde d_{j, k} \psi_{j, k}(x) +
\sum_{j = J_1}^{J - 1} \sum_{k = 0}^{2^j - 1}
\tilde d_{j, k} \psi_{j, k}(x).
\]

\smallskip

\noindent
By Lemma 3.4 we have for all $\omega' \in B_{\delta}$:

\smallskip 

\noindent (0) \ \ \ \  
$|c_{J_0 , k}^{(y)} - c_{J_0 , k}^{(f)}| = |c_{J_0 , k}^{(e(\omega'))}| \ \leq \ $
$\lambda_{n,\delta}$

\smallskip 

\noindent By Lemma \ref{tilda} we have for all $\omega' \in B_{\delta}$:

\smallskip

\noindent (1) \ \ \ \ 
$|\tilde d_{j, k} - d_{j, k}^{(f)}| \ \leq \ |d_{j,k}^{(f)}| \ \leq \ $
$C_{\varphi} M \cdot 2^{-j(\frac{1}{2} + \alpha)}$ \ \ \ \ \ \ \ \ \
for $J_1 \leq j < J$,\ $0 \leq k < 2^j$

\smallskip 

\noindent (2) \ \ \ \ 
$|\tilde d_{j, k} - d_{j, k}^{(f)}| \ \leq \ 2\lambda_{n,\delta} $
\ \ \ \ \ \ \ \ \ \ \ \ \ \ for $J_0 \leq j < J_1$, \ $0 \leq k < 2^j$.

\medskip

Let us first deal with the case of Haar wavelets (when $\alpha \leq 1$).
For a given $j$, the supports of different Haar wavelets do not overlap.
Therefore, for all $x \in \ ]0,1]$ there exist $K_1$ and $K(j)$ such that  

\smallskip
 
\noindent 
$|\tilde f(x) - \tilde y(x)| \ \leq \ $

$|c_{J_0 , K_1}^{(y)} - c_{J_0,K_1}^{(f)}| \cdot 2^{J_0/2} \ + \ $
$\sum_{j=J_0}^{J_1-1} |\tilde d_{j,K(j)} - d_{j,K(j)}^{(f)}| \cdot 2^{j/2} \
+ \ $
$\sum_{j=J_1}^{J-1} |\tilde d_{j,K(j)} - d_{j,K(j)}^{(f)}| \cdot 2^{j/2}$

\smallskip

\noindent This and (0), (1), (2) imply for all $x \in \ ]0,1]$:

\smallskip

$|\tilde f(x) - \tilde y(x)| \ \leq \ $
$C_1 \cdot \left( \frac{\log n}{n}\right) ^{\frac{\alpha}{1+ 2 \alpha}} + $
$C_2 \cdot \left( \frac{\log n}{n}\right)^{\frac{\alpha}{1+ 2 \alpha}}  + $
$C_3 \cdot \left( \frac{\log n}{n}\right)^{\frac{\alpha}{1+ 2 \alpha}}$ 

\smallskip

$= (c'_1 + c'_2 \sqrt{1+\delta}) \cdot$
$ \left( \frac{\log n}{n}\right)^{\frac{\alpha}{1+ 2 \alpha}}$

\smallskip

\noindent
Letting $x = \frac{i}{n}$ ($1 \leq i \leq n$) we obtain for all 
$\omega' \in B_{\delta}$:

\smallskip

$|f_i - \tilde y_i(\omega')| = |\tilde f(\frac{i}{n}) - \tilde
y(\frac{i}{n})|$
$ \ \leq \ $
$(c'_1 + c'_2 \sqrt{1+\delta}) \cdot $
$\left( \frac{\log n}{n}\right)^{\frac{\alpha}{1+ 2 \alpha}} $

\smallskip

\noindent
In the Haar case the theorem follows from this and the fact that 
$P(B_{\delta}) > 1 - \frac{9}{n^{1+ \delta}}$ (when $n \geq n_0$).

\bigskip

For wavelets on the interval (when $\alpha > 1$, 
and the number of vanishing moments is $N = \lceil \alpha \rceil$), 
there are never more than $2N$ wavelets that overlap (for a given $j$).
Indeed, in the above sums we have for each $j$ and each $x$: \ 
$0 \leq 2^j x - k \leq 2 N -1$. (Other values of $k$ would place the 
argument $2^j x - k$ of the wavelet functions outside of the support and
would ence only produce zero-terms in the sums.) Hence $k$ only needs to
range from \ $\lceil 2^j x \rceil - 2N +1$ \ through \ 
$\lceil 2^j x \rceil$, which corresponds to $2N$ values of $k$.   
 
Hence, the same calculation as for Haar wavelets applies, except that 
the constants $C_1$, $C_2$, $C_3$, $c'_1$, $c'_2$ need to be multiplied 
by $2N$. \ \ \ \ $\Box$

%%%%%%%%%%%%%%%%%%%%%%%%%%%%%%%%%%%%%%%%%%%%%%%%%%%%%%%%%%%%%%%%%%%%
\bigskip

\bigskip

\bigskip

\noindent {\large \bf  Appendix }

\bigskip

\noindent {\bf Proof of Lemma 2.1 }

\medskip

\noindent
Properties (H1) and (H2) follow from a direct calculation based on the exact
formulas for the Haar wavelets $\varphi_{j,k}$ and $\psi_{j,k}$.

\smallskip

$c^{(e)}_{j,k} = \int_0^1 \overline{e}(x) \varphi_{j,k}(x) dx =
2^{j/2} \int_{k2^{-j}}^{(k+1)2^{-j}} \overline{e}(x) dx =$

\smallskip

$\sum_{i=k2^{J-j}}^{(k+1)2^{J-j} -1} e_{i+1} 2^{-J} =
2^{-J+j/2} \sum_{i=0}^{2^{J-j}-1} e_{i + 1 + k 2^{J-j}}$.

\smallskip

\noindent The calculation for (H2) is similar.
The same calculation as for (H2) will give for $\overline{f}$: \

\smallskip

$d_{j,k}^{(f)} = 2^{-J-1+j/2}$
$\sum_{i=0}^{2^{J-j-1}} (f(i+1+ k 2^{J-j}) - f(i+1+ (k+ \frac{1}{2})
2^{J-j}))$. 

\smallskip

\noindent  
Then we use the H\"older condition \ \
$|f(i + 1 + k 2^{J-j}) - f(i + 1 + (k + \frac{1}{2}) 2^{J-j})| \ \leq \ $
$M \, (\frac{1}{2} 2^{J-j})^{\alpha}$.  \ \ \ \ \ \ \ \   $\Box$

\bigskip

\bigskip

\newpage
%%%%%%%%%%
\noindent {\bf Proof of Lemma 2.2}

\bigskip

\noindent Property (D1) follows from a direct calculation:

\smallskip

$c^{(e)}_{j,k} = \int_0^1 \overline{e}(x) \varphi_{j,k}(x) dx = $
$\sum_{i=0}^{n-1} e_i \int_{\frac{i}{n}}^{\frac{i+1}{n}} \varphi_{jk}(x) dx$

\smallskip

\noindent where we denote the functions $\varphi^{\mathrm{left}}_{jk}$ by
$\varphi_{j, 2^j-2N+k}$,  and $\varphi^{\mathrm{right}}_{jk}$ by
$\varphi_{j, 2^j-N+k}$.

\noindent
For the $\varphi_{jk}$ ``in the middle'' of the interval  we have

\smallskip

$\int_{\frac{i}{n}}^{\frac{i+1}{n}} \varphi_{jk}(x) dx$
$= 2^{j/2} \int_{i2^{-J+j}-k}^{(i+1)2^{-J+j}-k} \varphi(t) \, 2^{-j} dt$
$= 2^{j/2} 2^{-J} \alpha_{ijk}$

\smallskip

\noindent
by the the mean-value theorem, for some numbers $\alpha_{ijk}$ with \
$|\alpha_{ijk}| \ \leq \ \sup_{[0,1]} |\varphi|$.

\noindent
For the $\varphi_{j, 2^j-2N+k}$ ``at the left end'' of the interval,

\smallskip

$\int_{\frac{i}{n}}^{\frac{i+1}{n}} \varphi^{\mathrm{left}}_{jk}(x) dx$
$= \int_{\frac{i}{n}}^{\frac{i+1}{n}} \sum_{s=0}^{2N-1} (-s)^k \,
\varphi(2^jx+s) dx$
$= \sum_{s=0}^{2N-1} (-s)^k \int_{i2^{-J+j}+s}^{(i+1)2^{-J+j}+s} $
$\varphi(y) \, 2^{-j} dy$

\smallskip

$= \sum_{s=0}^{2N-1} (-s)^k 2^{-j} 2^{-J+j} \gamma_{ijs}$

\smallskip

\noindent
by the mean-value theorem, for some numbers $\gamma_{ijs}$ with
$|\gamma_{ijs}| \leq \sup_{[0,1]} |\varphi|$. By taking

\smallskip

$\alpha_{ijk} = 2^{-j/2} \sum_{s=0}^{2N-1} (-s)^k \gamma_{ijs}$

\smallskip

\noindent we obtain (D1).
At the left end, $k \leq N$, so \
$|\alpha_{ijk}| \leq  2N (2N-1)^N \cdot  \sup |\varphi|$.

\medskip

The scaling functions ``at the right end'' of the interval are handled in a
similar way.
The calculation for (D2) is similar.
(D3) follows from the wavelet characterization of H\"older classes
(\cite{Daubechies}, page 299, and \cite{Meyer}).  \ \ \ \ \ \ \ \   $\Box$

%%%%%%%%%%%%%%%%%%%%%%%%%%%%%%%%%%%%%%%%%%%%%%%%%%%%%%%%
% Bibliography

\end{document}